\documentclass[12pt]{article}
\usepackage{a4wide}
\usepackage{amstext,amssymb,amsbsy,array,graphicx,amsmath}

\def\bfF{{\boldsymbol F}}

\newcommand{\tnorm}[1]{\vert\hspace{-0.3mm}\Vert#1\Vert\hspace{-0.3mm}\vert}
\newcommand{\tn}{\vert\hspace{-0.3mm}\Vert}
\newcommand{\bfv}{{\boldsymbol v}}
\newcommand{\bfpsi}{{\boldsymbol\psi}}

\newcommand{\bfpi}{{\boldsymbol \pi}}

\newcommand{\bfsig}{{\boldsymbol\sigma}}

\newcommand{\bftheta}{{\boldsymbol \theta}}
\newcommand{\bfzet}{{\boldsymbol \zeta}}

\newcommand{\bftwei}{{\boldsymbol \vartheta}}

\newcommand{\bfeta}{{\boldsymbol\eta}}

\newcommand{\pit}{{\boldsymbol P}}
\newcommand{\piz}{{\boldsymbol Q}}
\newcommand{\bfR}{{\boldsymbol R}}

\newcommand{\bfI}{{\boldsymbol 1}}

\newcommand{\bbR}{{\mathbb{R}}}

\newcommand{\mcE}{\mathcal{E}}
\newcommand{\mcT}{\mathcal{T}}
\newcommand{\mcK}{\mathcal{K}}

\newcommand{\bfV}{{\boldsymbol V}}
\newcommand{\bfn}{{\boldsymbol n}}

\newcommand{\bfeps}{{\boldsymbol\varepsilon}}

\newcommand{\bfx}{{\boldsymbol x}}

\newcommand{\norm}[2]{\|#1\|_{#2}}

\newtheorem{theorem}{Theorem}[section]
\newtheorem{thm}{Theorem}[section]
\newtheorem{lemma}[theorem]{Lemma}
\newtheorem{lem}[theorem]{Lemma}
\newtheorem{remark}[theorem]{Remark}

\newenvironment{pf}{\noindent \newline {\bf Proof.}}{\hfill \mbox{\fbox{} } \newline}

\title{\bf Locking Free Quadrilateral Continuous/Discontinuous Finite Element 
Methods for the Reissner--Mindlin Plate} 
\author{Peter Hansbo and Mats G. Larson}
\begin{document}
\maketitle

\begin{abstract}
{We develop a finite element method with continuous displacements and 
discontinuous rotations for the Mindlin-Reissner plate model on 
quadrilateral elements. To avoid shear locking, the rotations 
must have the same polynomial degree in the parametric reference plane as 
the parametric derivatives of the displacements, and obey the same transformation 
law to the physical plane as the gradient of displacements. We prove optimal 
convergence, uniformly in the plate thickness, and provide numerical results 
that confirm our estimates.}
\end{abstract}

\section{Introduction}

\paragraph{The Reissner-Mindlin Plate Model and Shear Locking.} 
The Reissner-Mindlin equations is a model of the displacement of a moderately
thick plate under transversal load. The unknows are the normal displacement field 
$u$ and the rotation field $\bftheta$ of a normal fiber. The difficulty with this 
model, from a numerical point of view, is the matching of the approximating 
spaces for $\bftheta$ and $u$. As the thickness $t\rightarrow 0$, the difference 
$\nabla u-\bftheta$ must tend to zero, which, for naive choices of spaces, 
leads to a deterioration of the approximation known as locking or in this case 
shear locking since the difficulty emanates from the term involving the shear 
energy. The situation is particularly difficult if we wish to use low order 
approximations. 

\paragraph{Earlier Work.}
There are basically three different approaches to solve this problem. 
Perhaps the most common approach has been to use a projection to relax 
the equation, which essentially corresponds to a mixed formulation where an 
additional variable, often the shear vector proportional to $(\nabla u - \bftheta)/t^2$, is introduced.  
For instance, the MITC element family of Bathe and co-workers \cite{BaDv85} are 
based on this approach. For quadrilaterals, this type of approach has been
used and analyzed in \cite{ArBo02b, DuLi04, DuHe03, MiSh05,Ye00}.

Another approach is to use a stabilized mixed formulation, see 
Chapelle and Stenberg, \cite{ChSt99a, ChSt99b}. 

Finally, a third approach is to use finite element spaces that 
are rich enough to satisfy the shear constraint exactly while maintaining optimal 
approximation properties. This approach was first proposed by Hansbo and Larson 
\cite{HaLa03}, where continuous piecewise quadratics for the displacements and 
discontinuous piecewise linears for the rotations in a discontinuous Galerkin 
formulation. Further developments, still using simplicial elements, were given by 
Arnold et al. \cite{ArBr07}, Heintz et al. \cite{HaHe11}, and B\"osing et al. 
\cite{BoMa10}. When the thickness of the plate tends to zero we obtain the Kirchhoff 
plate and our scheme can be 
seen as a version of the method proposed in \cite{EnGa02}, see \cite{HaHe11}. In 
this context we also mention the fully discontinuous Galerkin method developed 
in \cite{HaLa02} and the parametric continuous/discontinuous Galerkin method 
\cite{BrSe13} for the Kirchhoff plate.

\paragraph{New Contributions.}
In this paper we extend the method of \cite{HaLa03} to quadrilateral elements. We 
show that, with the proper definition of the finite element space for the rotations, 
we can satisfy the equation $\nabla u - \bftheta = 0$ exactly while maintaining optimal 
approximation properties and thus together 
with stability we obtain optimal a priori error estimates uniformly in the 
thickness parameter. Using continuous tensor product quadratics for the displacements 
the suitable space for the rotations consists of discontinuous parametric vector 
polynomials that are also mapped in the same way as the gradient of elements. The 
mapping is the rotated, or covariant, Piola mapping that preserves tangent 
traces, and naturally appears in the context of curl conforming elements, see 
\cite{FaGa11}. We could also use the smaller subspace of tangentially continuous functions 
for the rotations instead of the full discontinuous space. The interpolation error 
estimates on quadrilaterals are based on the observation that tensor product polynomials 
mapped with a bilinear map contain complete polynomials, cf. \cite{ArBo02,ArBo05}, and thus 
the estimates follows from the Bramble--Hilbert lemma and scaling. We also show that the condition 
required to avoid locking also implies that complete linear polynomials are contained 
in the space for the rotations. Our analysis is remarkably simple and avoids difficulties 
caused by the mixed formulations and basically rely on proper construction of the discrete 
spaces and approximation properties that takes advantage of the stability of the underlying 
continuous problem.

We remark that the idea of using covariant maps to obtain suitable approximations of the rotations 
has also recently been used in the context of isogeometric approximations
by Beir{\~a}o~da Veiga, Buffa, Lovadina, Martinelli, and Sangalli \cite{BeBuLoMa12}.

\paragraph*{Outline.} In Section 2 we formulate the Reissner--Mindlin model on weak form, in Section 3 we introduce the quadrilateral finite element spaces 
and formulate the finite element method, in Section 4 we derive approximation properties and a priori 
error estimates, and in Section 5 we present numerical results illustrating the theoretical 
results.

\section{The Reissner-Mindlin Plate Model}
\subsection{Energy Functional}
Consider a plate with thickness $t$ occupying a convex polygonal domain $\Omega$ in ${\bbR}^2$, 
which is clamped at the boundary $\partial \Omega$. The Reissner-Mindlin plate model can be derived 
from minimization of the sum of the bending energy, the shear energy, and the potential of the surface 
load
\begin{equation}\label{min}
{\mathcal{E}_\text{RM}}(u,\bftheta) :=
\frac{1}{2}a(\bftheta,\bftheta)+\frac{\kappa}{2\,
t^2}\int_{\Omega}\vert \nabla u-\bftheta\vert^2\, d\Omega
-\int_{\Omega} g\, u\, d\Omega 
\end{equation}
Here $u$ is the transverse displacement, $\bftheta$ is the rotation of the median 
surface, $t$ is the thickness, $t^3\, g$ is the transverse surface load, and the 
bending energy $a(\cdot,\cdot)$ is defined by
\begin{equation}
a(\bftheta,\bftwei) := \int_{\Omega}\Bigl(
2\mu\bfeps(\bftheta):\bfeps(\bftwei)+\lambda \nabla\cdot\bftheta\,
\nabla\cdot\bftwei\Bigr)\, d\Omega
\end{equation}
where $\bfeps$ is the curvature tensor
\begin{equation}
\bfeps(\bftheta) : = \left[\begin{array}{>{\displaystyle}c>{\displaystyle}c}
\frac{\partial\theta_x}{\partial x}& \frac12\left(\frac{\partial\theta_x}{\partial y}+\frac{\partial\theta_y}{\partial x}\right)\\[4mm]
\frac12\left(\frac{\partial\theta_x}{\partial y}+\frac{\partial\theta_y}{\partial x}\right) &\frac{\partial\theta_y}{\partial y}\end{array}\right] .
\end{equation}
The material parameters are
given by the relations $\kappa = E\,k/(2(1+\nu))$, $\mu := E/(24
(1+\nu))$, and $\lambda := \nu E /(12 (1-\nu^2))$, where $E$ and
$\nu$ are the Young's modulus and Poisson's ratio, respectively,
and 
$k$ is a shear correction factor. We shall alternatively write the 
bending energy product as
\begin{equation}
a(\bftheta,\bftwei) = \int_\Omega \bfsig(\bftheta):\bfeps(\bftwei)\, d\Omega
\end{equation}
where $\bfsig(\bftheta) := 2\mu\bfeps(\bftheta) + \lambda \nabla\cdot\bftheta\, \bfI$
is the moment tensor.


\subsection{Weak Form}
The transverse displacement and rotation vector are
solutions to the following variational problem: find $(u, \bftheta ) 
\in  H_0^1(\Omega) \times [H_0^1(\Omega)]^2$ such that
\begin{equation}\label{exactsol}
a(\bftheta, \bftwei) 
+ \frac{\kappa}{t^2}(\nabla u-\bftheta, \nabla v-\bftwei ) = (g, v), 
\quad \forall (v,\bftwei) \in H_0^1(\Omega) \times [H_0^1(\Omega)]^2
\end{equation}
where $(\cdot,\cdot)$ denotes the $L^2(\Omega)$ inner product, 
$H^k(\Omega)$ are the usual Sobolev spaces, and the functions in 
$H^1_0(\Omega)$ have zero trace on the boundary $\partial \Omega$.

\section{The Finite Element Method}
\subsection{The Quadrilateral Mesh}
Next, let $\mcK^h = \{ K \}$ be a family of quasiuniform partitions of $\Omega$ into convex 
quadrilaterals $K = \bfF_K(\widehat{K})$ with mesh parameter $h$ such that $ c h \leq h_K \leq C h$, 
where $h_K = \text{diam}(K)$,  for all $K \in \mcK^h$. We also assume that $\mcK^h$ is a 
shape regular partition in the sense that $h_K/\rho_K \leq C$ for all $K \in \mcK^h$, where 
$\rho_K$ is the smallest diameter of the largest inscribed circle in any of the four subtriangles 
obtained by inserting a diagonal between two opposite corners in $K$.

\subsection{Parametric Elements for Displacements and Rotations}

\newcommand{\hatT}{\widehat{T}}
\newcommand{\hatP}{\widehat{P}}
\newcommand{\hatSigma}{\widehat{\Sigma}}
\newcommand{\hatx}{\widehat{x}}
In order to define our finite element spaces we begin with a continuous 
parametric finite element space $V_{\rm D}^h$ for the displacement $u$ 
and then we determine a space $\bfV_{\rm R}^h$ of discontinuous piecewise 
parametric functions for the rotations $\bftheta$ such that 
\begin{equation}
\nabla V_{\rm D}^h \subseteq \bfV_{\rm R}^h
\end{equation}
in order to be able to satisfy the equation $\bftheta - \nabla u =  0$ 
exactly, when the thickness tends to zero. Using this inclusion we 
identify the proper space for the rotations. For clarity, we restrict 
the presentation to quadratic tensor product approximation of the 
displacements. extension to higher order elements follow directly.

\subsubsection{Displacements}

Let $\widehat{K}$ be the reference unit square and $Q_{k,l}(\widehat{K})$ 
the space of tensor product polynomials of order $k$ and $l$ in each variable, 
more precisely
\begin{align}
Q_{k,l}(\widehat{K}) &= \text{span} \{ \widehat{x}^\alpha \widehat{y}^\beta\; : \; 0 \leq \alpha \leq k, 0\leq  \beta \leq l \}
\end{align}
and $Q_{k}(\widehat{K}) = Q_{k,k}(\widehat{K})$. For each $K\in \mcK^h$ 
let $\bfF_K:\widehat{K} \rightarrow {\bf R}^2$ 
be the bilinear, i.e., $\bfF_K \in [Q_1(\widehat{K})]^2$, mapping such that 
$K=\bfF_K(\widehat{K})$. 
We define the space of parametric tensor product polynomials on $K$ by
\begin{equation}
V_{{\rm D},K} = \{ p:K \rightarrow {\bf R}\, : \, p = \widehat{p} \circ \bfF_K^{-1} \}
\end{equation}
and the corresponding space on $\mcK^h$ of continuous piecewise parametric 
tensor product polynomials  
\begin{equation}
V^h_{\rm D} = \{v: \Omega \rightarrow \bfR : v |_K \in V_{{\rm D},K} \; \forall K \in \mcK^h, v \in C(\Omega) \}
\end{equation}

\subsubsection{Rotations}

Turning to the space for the rotations we recall that, since 
$p = \widehat{p} \circ \bfF_K^{-1}$, we have 
\begin{equation}
\nabla p (x) = DF^{-\rm{T}}_K \widehat{\nabla} \widehat{p}(\widehat{x} ) 
=DF^{-\rm{T}}_K \widehat{\nabla} \widehat{p}(\bfF^{-1}_K(x)) 
\end{equation}
where $\widehat{\nabla}$ is the gradient in the reference coordinates. 
Introducing the rotated or covariant Piola mapping 
\begin{equation}
\bfR_K:\widehat{\bfV}_{\rm R} \ni \widehat{\bftheta} \mapsto D\bfF^{-{\rm T}}_K \widehat{\bftheta} \circ \bfF_K^{-1} \in {\bfV_{\rm R}}
\end{equation}
we have $\nabla p (x) = DF^{-\rm{T}}_K \widehat{\nabla} \widehat{p}(\widehat{x} ) 
= \bfR_K \widehat{\nabla} \widehat{p}(\widehat{x} )$.  We 
are thus led to defining the following space for the rotations 
\begin{equation}
\bfV_{{\rm R},K} = \bfR_K {\widehat{\bfV}}_{\rm R} 
\end{equation}
where ${\widehat{\bfV}}_{\rm R} $ is a space on the reference unit 
square $\widehat{K}$ that satisfies
\begin{equation}
\widehat{\nabla} {\widehat{V}}_{{\rm D},K}  = Q_{1,2}(\widehat{K}) \times Q_{2,1}(\widehat{K}) \subseteq  \widehat{\bfV}_{{\rm R},K}
\end{equation}
We finally define the space of discontinuous mapped parametric functions
\begin{equation}
\bfV^h_{\rm R} = \{\bftheta: \Omega \rightarrow {\bfR}^2  
: \bfv|_K \in \bfV_{{\rm R},K} \; \forall K \in \mcK^h \}
\end{equation}
\begin{remark}
We note that it is indeed also possible to chose a subspace of $\bfV^h_{\rm R}$ 
that consists of functions that have continuous tangential trace at each edge. This 
case is also covered by our analysis and basically only depends on the choice of interpolation 
operator on the reference element.
\end{remark}

\subsection{The Finite Element Method}

Let $\mcE^h =\{ E \}$ be the set of edges in the mesh $\mcK^h$. We 
split $\mcE^h$ into two disjoint subsets
\begin{equation}
  \mcE^h = \mcE^h_I \cup \mcE^h_B
\end{equation}
where $\mcE^h_I$ is the set of edges in the interior of $\Omega$ and
$\mcE^h_B$ is the set of edges on the boundary $\partial \Omega$. 
Further, with each edge we associate a fixed unit normal $\bfn$ 
such that for edges on the boundary $\bfn$ is the exterior unit 
normal. We denote the jump of a function $\bfv \in \bfV_{\rm R}^h$ at 
an edge $E$ by $ \left[\bfv \right] = \bfv^+-\bfv^-$ for 
$E \in \mcE^h_I$ and $\left[\bfv\right] = \bfv^+$ for $E \in \mcE^h_B$, 
and the average
$\langle  \bfv\rangle = (\bfv^+ + \bfv^-)/2$ for $E \in \mcE^h_I$
and $\langle \bfv\rangle = \bfv^+$ for $E \in \mcE^h_B$, where
$\bfv^{\pm}(\bfx) = \lim_{\epsilon\downarrow 0} \bfv(\bfx\mp
\epsilon\,\bfn)$ with $\bfx\in E$.

The method takes the form:
find $(u^h,\bftheta^h) \in V^h_{\rm D} \times \bfV^{h}_{\rm R}$ such that
\begin{equation}\label{niteq}
a_h(\bftheta^h,\bftwei)+ \frac{\kappa}{t^2}\left(\nabla u^h-\bftheta^h,\nabla
v-\bftwei\right)=(g,v) \quad \forall (v,\bftwei) \in V^h_{\rm D}\times \bfV^h_{\rm R}
\end{equation} 
Here the bilinear form $a_h(\cdot,\cdot)$ is defined by
\begin{equation} \label{nitsche_form}\begin{array}{>{\displaystyle}l>{\displaystyle}l}
a_h(\bftheta^h,\bftwei) &= \sum_{K \in \mcK^h}
(\bfsig(\bftheta^h), \bfeps( \bftwei) )_K \\ \nonumber &
\quad - \sum_{E \in \mcE^h_I \cup \mcE^h_B } ( \langle\bfn \cdot 
\bfsig(\bftheta^h) \rangle, [ \bftwei ])_E
     + (\langle \bfn \cdot \bfsig(\bftwei)\rangle, [\bftheta^h ] )_E \\ \nonumber
& \quad + (\mu + \lambda)\, \gamma \sum_{E \in \mcE^h_I \cup
\mcE^h_B}
         ( h_E^{-1} [ \bftheta^h ],[ \bftwei] )_E
\end{array}
\end{equation}
where $\gamma$ is a positive constant, $h_E$ is defined by
\begin{equation}\label{hedge}
h_E =
 \left( |T^+| + |T^-| \right) / (2 \, |E|) \quad \text{for $E = \partial T^+ \cap
\partial T^-$}
\end{equation}
with $|K|$ the area of $K$, on each edge, and
$(\cdot, \cdot)_\omega$ is the $L^2(\omega)$ inner product with $\omega = K, E$.


\section{A Priori Error Estimates}

The analysis presented here extends the analysis in Hansbo and Larson \cite{HaLa03} 
to parametric elements on quadrilaterals. For completeness we include the necessary 
results but refer to \cite{HaLa03} and \cite{HaHe11}, for further details.

\subsection{Stability and Continuity of the Discrete Bilinear Form}

Let the mesh dependent energy-like norm, associated with the 
bilinear form $a_h(\cdot,\cdot)$ be defined by
\begin{align}\label{energy_norm}
\tnorm{\bftwei}^2  = & \sum_{K \in \mcK^h} (\bfsig(\bftwei),
\bfeps(\bftwei))_K  +\sum_{E \in \mcE^h} \frac{1}{2 \mu + 2\lambda}\|
 h_E^{1/2}  \langle \bfn \cdot \bfsig(\bftwei)\rangle\|_E^2 
 \\ \nonumber
&  \qquad +\sum_{E \in \mcE^h}  (2 \mu + 2\lambda)\|
       h_E^{-1/2} \left[ \bftwei\right]\|^2_E
\end{align}
We summarize the standard properties in the following and then 
state Cea's lemma.
\begin{lem}\label{lemma:infsup} It holds:\\
\begin{itemize}
\item {\bf Continuity:} There is a constant $C$ such that
\begin{align}
&a_h(\bftheta,\bfeta ) + \frac{\kappa}{t^2} (\nabla u - \bftheta, \nabla v - \bftwei ) 
\\ \nonumber
&\qquad \leq 
C \left( \tn \bftheta \tn^2 + \frac{\kappa}{t^2} \| \nabla u - \bftheta\|^2\right)^{\frac{1}{2}}
\left( \tn \bfeta \tn^2 + \frac{\kappa}{t^2} \| \nabla v - \bftwei\|^2\right)^{\frac{1}{2}}
\end{align}
for all $(u,\bftheta), (v,\bftwei) \in H^1_0(\Omega) \times ( [H^1_0(\Omega)]^2 +  \bfV_{\rm R}^h )$.
\item {\bf Coercivity:} For $\gamma$ sufficiently large, there is a constant
$m$ such that
\begin{equation}\label{inter21}
m \left( \tn \bftheta \tn^2 + \frac{\kappa}{t^2} \| \nabla u - \bftheta\|^2 \right)
\leq a_h(\bftheta,\bftheta ) + \frac{\kappa}{t^2} (\nabla u - \bftheta, \nabla u - \bftheta )
\end{equation}
for all $(u,\bftheta) \in V^h_{\rm D} \times \bfV^h_{\rm R}$. 

\item{\bf Consistency:} With $(u,\bftheta)$ the exact solution defined by (\ref{exactsol}) 
and $(u^h,\bftheta^h)$ the discrete solution defined by (\ref{niteq}) we have 
\begin{equation}\label{galort}
a_h(\bftheta - \bftheta^h,\bftwei ) + \frac{\kappa}{t^2} (\nabla (u - u^h) 
- (\bftheta - \bftheta^h), \nabla v - \bftwei ) = 0
\end{equation} 
for all $(v,\bftwei) \in V^h_{\rm D}\times \bfV^h_{\rm R}$. 
\end{itemize}
\end{lem}
\begin{pf} The continuity estimate follows directly from Cauchy-Schwartz. The coercivity 
follows from coercivity for $a_h(\cdot,\cdot)$, which depend on the inverse inequality 
\begin{equation}
\frac{1}{2 \mu + 2 \lambda } \|h_E^{\frac{1}{2}} \bfsig(\bftwei)\|^2_{\partial K} 
\leq C (\bfsig(\bftwei), \bfeps(\bftwei) )_K \quad \forall \bftwei \in \bfV_{{\rm R},K}
\end{equation}
This inequality is established by mapping to the reference element and using finite 
dimensionality and then mapping back in the same way as for affine elements. The 
consistency  follows by using Green's formula.
\end{pf}
\begin{lem}\label{lem:cea} (Cea's Lemma) With $(u,\bftheta)$ the exact solution defined by (\ref{exactsol}) 
and $(u^h,\bftheta^h)$ the discrete solution defined by (\ref{niteq}) we have 
\begin{align}\nonumber
&\tn \bftheta - \bftheta^h \tn^2 + \frac{\kappa}{t^2} \| \nabla (u - u^h) - (\bftheta - \bftheta^h)\|^2
\\ \label{cea}
&\qquad \leq C \left( \tn \bftheta - \bftwei \tn^2 + \frac{\kappa}{t^2} \| \nabla (u - v) - (\bftheta - \bftwei)\|^2 \right)
\end{align} 
for all $(v,\bftwei) \in V^h_{\rm D} \times \bfV^h_{\rm R}$. 
\end{lem}
\begin{pf} This estimate follows by first splitting the error $(u,\bftheta) - (u^h,\bftheta^h) 
= (u,\bftheta) - (v,\bftwei) + (v, \bftwei) - (u^h,\bftheta^h)$ and then using coercivity 
followed by consistency and finally continuity for the second term. 
\end{pf}

\subsection{Interpolation}

\subsubsection{Parametric Lagrange Interpolation}

Let $\pi_{{\rm D},K}: H^2(K) \rightarrow V_{{\rm D},K}$  
be the parametric Lagrange interpolant 
defined by 
\begin{equation}
\pi_{{\rm D},K} v  
= (\widehat{\pi}_{\rm D} \widehat{v}) \circ \bfF_K^{-1}, 
\quad \forall K \in \mcK^h 
\end{equation}
where $\widehat{\pi}_{\rm D}$  is the usual nodal 
Lagrange interpolant on the reference element. We then have the following interpolation 
error estimate. The short proof, essentially following \cite{ArBo02, ArBo05}, is included 
and will be reused when we consider approximation 
properties for the rotations. 

\begin{lemma}\label{lemma:interpolD} The following estimate holds
\begin{equation}\label{interpolscalar}
\|u - \pi_{{\rm D},K} u \|_{K,m} \leq C h^{k - m} \| u \|_{K,k}, \quad 0 \leq m \leq k \leq 3   
\end{equation}
\end{lemma}

\begin{pf}  We first prove the estimate in the case $h_K =1$ and then obtain 
the general estimate using scaling. Let $P_k(K)$ denote the space of polynomials 
of order less or equal to $k$. The key observation is that $P_2(K) 
\subseteq V_{{\rm D},K}$, since $\widehat{V}_{\rm D}$ contains $Q_2(\widehat{K})$.  This fact 
follows directly by observing that $ P_2(K) \circ \bfF_K \subset Q_2(\widehat{K})$ 
and thus $P_2(K)\subseteq Q_2(\widehat{K})\circ \bfF_K^{-1}$. We also note that 
$\pi_{{\rm D},K} p =  p$ for all $p \in P_2(K)$ since $\widehat{\pi}_{\rm D}$ is a projection 
onto $Q_2(\widehat{K})$. We thus have
\begin{align}
\|u - \pi_{{\rm D},K} u \|_{K,m} &\leq \inf_{p \in P_2(K)} \| (I - \pi_{{\rm D},K}) ( u - p ) \|_{K,m} 
\\
&\leq \| I - \pi_{{\rm D},K} \|_{\mathcal{L}(H^k(K),H^m(K))} \inf_{p \in P_2(K)}  \|  u - p \|_{K,k}
\\
&\leq C \| I - \pi_{{\rm D},K} \|_{\mathcal{L}(H^k(K),H^m(K))}  |  u |_{K,k}
\end{align}
where we used the Bramble-Hilbert lemma in the last inequality. Using shape regularity we have 
$\| I - \pi_{{\rm D},K} \|_{\mathcal{L}(H^k(K),H^m(K))} \leq C$ since 
$\| I - \widehat{\pi}_{\rm D} \|_{\mathcal{L}(H^k(\widehat{K}),H^m(\widehat{K}))} \leq C$
and thus the result follows in the case 
$h_K=1$. Finally, using the dilation $\bfx \mapsto h_K^{-1} \bfx$ we map an arbitrary 
quadrilateral $K$ to a quadrilateral $\tilde{K}$ with $h_{\tilde{K}}=1$. We then have
\begin{align}
\|u - \pi_{{\rm D},K} u \|_{K,m} &= h_K^{1-m} \| u - \pi_{{\rm D},K} u \|_{\tilde{K},m}
\\
&\leq C h_K^{1-m} | u |_{\tilde{K},k}
\\
&= C h^{k-m} | u |_{K,k}
\end{align}
which completes the proof. 
\end{pf}

\subsubsection{Interpolation for Rotations}
Let $\bfpi_{\rm R}: [H^2(\Omega)]^2 \rightarrow \bfV_{\rm R}^h$ be defined by
\begin{equation}
\bfpi_{{\rm R},K} \bfeta = \bfR_K  \widehat{\bfpi} \bfR_K^{-1}\bfeta 
= D\bfF_K^{-\rm{T}} ( \widehat{\bfpi}(D\bfF_K^{\rm{T}} \bfeta\circ \bfF_K ))\circ \bfF_K^{-1}
\end{equation}
Then we first have the following interpolation error estimate.
\begin{lem}\label{lemma:interpolR} The following estimate holds
\begin{equation}\label{interpolvector}
\|\bfeta  - \bfpi_{{\rm R},K} \bfeta \|_{K,m} \leq C h^{2 - m} \|\bfeta \|_{K,2} \quad m =0,1
\end{equation}
\end{lem}
\begin{pf} Starting from the fact that $P_2(K) \subset V_{{\rm D},K}$ we have 
\begin{equation}
[P_1(K)]^2 \subseteq \nabla P_2(K) \subseteq \nabla V_{{\rm D},K} 
=  \bfR_K \widehat{\nabla} \widehat{V}_{\rm D}
\subseteq \bfR_K \widehat{\bfV}_{\rm R} = \bfV_{{\rm R},K}
\end{equation}
where we used the inclusion $\widehat{\nabla}\widehat{V}_{\rm D} \subseteq \widehat{\bfV}_{\rm R}$, 
and thus we conclude that
 \begin{equation}
[P_1(K)]^2 \subseteq \bfV_{{\rm R},K}
\end{equation}
We may now prove the estimate using the same technique as in Lemma 
\ref{lemma:interpolD}.
\end{pf}

\begin{remark} We note that our results can be directly generalized to the following 
situation: If
\begin{equation}
Q_k(\widehat{K}) \subseteq \widehat{V}_{\rm D}, 
\qquad 
\widehat{\nabla} \widehat{V}_{\rm D} \subseteq \widehat{\bfV}_{\rm R}
\end{equation}
then the spaces
\begin{equation}
V_{\rm D} = \widehat{V}_{\rm D} \circ \bfF^{-1}_K, \qquad \bfV_{\rm R} = \bfR_K \widehat{\bfV}_{\rm R}
\end{equation}
have optimal interpolation properties. The condition for avoiding locking thus 
also implies optimal interpolation properties for the rotations.
\end{remark}

\subsection{Energy Norm A Priori Error Estimate}

\subsubsection{The Shear Stress}

We define the scaled shear stress $\bfzet$ and its discrete counterpart
 $\bfzet^h$, as follows
\begin{equation}\label{zeta}
\bfzet := \kappa^{1/2}(\nabla u -\bftheta)/t^2 
\quad\text{and}\quad
\bfzet^h := \kappa^{1/2}(\nabla u^h -\bftheta^h)/t^2
\end{equation}
Note that $\bfzet^h \in \bfV_{\rm R}^h$ due to the inclusion $\nabla V_{\rm D}^h \subset \bfV_{\rm R}^h$.
\subsubsection{A Stability Estimate}

Splitting the Reissner-Mindlin displacement $u = u_0 + u_r$, with 
$u_0$ the Kirchhoff solution obtained in the limit case $t \rightarrow 0$ 
and $u_r$ the difference between the solutions, we have the following stability 
estimate.
\begin{lem}\label{chapelle}
Assume that $\Omega$ is convex and $g\in L_2(\Omega)$. Then it holds
\begin{equation}
\norm{u_0}{H^3(\Omega)}+\frac{1}{t}\norm{u_r}{H^2(\Omega)}+\norm{\bftheta}{H^2(\Omega)}
+t\norm{\bfzet}{H^1(\Omega)}\leq
 C\Bigl(\norm{g}{H^{-1}(\Omega)}+t\norm{g}{L_2(\Omega)}\Bigr)
\end{equation}
\end{lem}
\begin{pf} See \cite{ArFa89} and \cite{ChSt98}. \end{pf}

\subsubsection{Approximation}

In order to make use of the stability result (\ref{chapelle}) in Cea's Lemma \ref{lem:cea} we 
introduce the operators
$\pit :  [H^2(\Omega) ]^2\rightarrow \bfV^h_{\rm R}$ and $\piz : [H^2(\Omega) ]^2\rightarrow \bfV^h_{\rm R}$ 
defined by
\begin{equation}
\pit\bftheta := \nabla \pi_{\rm D} u_0 -\bfpi_{\rm R} \nabla u_0 + \bfpi_{\rm R} \bftheta
\end{equation}
and
\begin{equation}
\piz\bfzet :=\kappa^{1/2}t^{-2}\left( \nabla\pi_{\rm D} u_r -\bfpi_{\rm R} \nabla u_r)\right) + \bfpi_{\rm R} \bfzet
\end{equation}
We then have the following two lemmas.

\begin{lem}\label{lemma:error} With $(u,\bftheta)$ the exact solution defined by 
(\ref{exactsol}), $(u^h,\bftheta^h)$ the discrete solution defined by (\ref{niteq}), and 
the continuous and discrete shear stress, $\bfzet$ and $\bfzet^h$, defined by 
(\ref{zeta}), we have the estimate 
\begin{align}
\tn \bftheta - \bftheta^h \tn^2 + t^2 \| \bfzet - \bfzet^h \|^2 
&\leq C \left( \tn \bftheta - \pit \bftheta \tn^2 + t^2 \| \bfzet - \piz \bfzet\|^2 \right)
\end{align} 
\end{lem}
\begin{pf}
Setting $v = \pi_{\rm D} u$ and $\bftwei = \pit \bftheta$ in (\ref{cea}) we have
\begin{align}\nonumber
&\tn \bftheta - \bftheta^h \tn^2 + \frac{\kappa}{t^2} \| \nabla (u - u^h) - (\bftheta - \bftheta^h)\|^2
\\ 
&\qquad \leq C \left( \tn \bftheta - \pit\bftheta \tn^2 
+ t^2 \| \kappa^{1/2}t^{-2} (\nabla (u - \pi_{\rm D} u ) - (\bftheta - \pit \bftheta) )\|^2 \right)
\end{align} 
and
\begin{align}
& \nabla (u - \pi_{\rm D} u ) - (\bftheta - \pit \bftheta)  \nonumber
\\
&\qquad = \nabla (u - \pi_{\rm D} u ) 
- (\bftheta - (\nabla \pi_{\rm D} u_0 - \bfpi_{\rm R} \nabla u_0 ) - \bfpi_{\rm R} \bftheta ) )
\\
&\qquad = (\nabla u - \bftheta ) - \bfpi_{\rm R} (\nabla u - \bftheta )
\\ \nonumber
&\qquad \qquad + (\nabla \pi_{\rm D} u_0 - \bfpi_{\rm R} \nabla u_0 ) +  \bfpi_{\rm R} \nabla u - \nabla \pi_{\rm D} u
\\
&\qquad = (\nabla u - \bftheta ) - \bfpi_{\rm R} (\nabla u - \bftheta )
\\ \nonumber
&\qquad \qquad +  \bfpi_{\rm R} \nabla u_r - \nabla \pi_{\rm D} u_r 
\\ 
&=\kappa^{-1/2} t^2 (\bfzet - \piz \bfzet )
\end{align}
which concludes the proof.
\end{pf} 

\begin{lem}\label{lemma:app}
We have the following estimate
\begin{align}
&\tnorm{\bftheta-\pit\bftheta}+  t \norm{\bfzet -\piz
\bfzet}{L_2(\Omega)} \nonumber
\\
&\qquad \leq C h \Bigl( \norm{\bftheta}{H^{2}(\Omega)}+\norm{u_0}{H^{3}(\Omega)}
 + t^{-1}\norm{u_r}{H^{2}(\Omega)} + t \norm{\bfzet}{H^{1}(\Omega)}\Bigr)
\end{align}
\end{lem}
\begin{pf}
To estimate $\tnorm{\bftheta-\pit\bftheta}$ we employ the trace inequality 
\begin{equation}\label{trace}
{h_K^{-1}}\norm{\bftwei}{L_2(\partial K)}^2
\leq C \left(h_K^{-2}\norm{\bftwei}{L_2(K)}^2+\norm{\bftwei}{H^1(K)}^2\right)
\quad \forall \bftwei \in [H^2(K)]^2
\end{equation}
to get the estimate
\begin{equation}
\tn \bftheta - \pit \bftheta \tn^2 
\leq \sum_{K \in \mcK^h} \sum_{l=0}^2 h_K^{2(l-1)} \| \bftheta - \pit \bftheta \|^2_{K,l}
\end{equation}
By the definition of $\pit$ and the triangle inequality we have 
\begin{align}
\norm{\bftheta-\pit\bftheta}{K,l}&\leq\norm{\bftheta-\bfpi_{\rm R} \bftheta}{K,l}+
\norm{\nabla u_0-\nabla \pi_{\rm D} u_0}{K,l}+\norm{\nabla u_0-\bfpi_{\rm R} \nabla
u_0}{K,l} 
\\
&\leq C h_K^{2-l} \Bigl(\norm{\bftheta}{H^{2}(K)} +\norm{u_0}{H^{3}(K)}\Bigr)
\end{align}
where we used the interpolation estimates (\ref{interpolscalar}) 
and (\ref{interpolvector}) in the second inequality. 
Thus we obtain the estimate
\begin{equation}
\tn \bftheta - \pit \bftheta \tn \leq C h \Bigl(\norm{\bftheta}{H^{2}(\Omega)} 
+ \norm{u_0}{H^{3}(\Omega)}\Bigr)
\end{equation}

Next we estimate the second term $\norm{\bfzet-\piz\bfzet}{L_2(\Omega)}$ using the definition of 
$\piz$, the triangle inequality, and the interpolation estimates (\ref{interpolscalar}) 
and (\ref{interpolvector}), as follows
\begin{align*}
t \norm{\bfzet-\piz\bfzet}{L_2(\Omega)}   \leq {} &
t \norm{\bfzet-\bfpi_{\rm R} \bfzet}{L_2(\Omega)}+\frac{\kappa^{1/2}}{t}\norm{\nabla u_r-  \nabla \pi_{\rm D}
u_r}{L_2(\Omega)}
\\ {} & \quad + \frac{\kappa^{1/2}}{t}\norm{\nabla u_r-\bfpi_{\rm R} \nabla u_r}{L_2(\Omega)}
\\
   \leq {} & C h \Bigl(t^{-1}\norm{u_r}{H^{2}(\Omega)}+ t \norm{\bfzet}{H^{1}(\Omega)}\Bigr)
\end{align*}
which completes the proof of the lemma.
\end{pf}

\subsubsection{Error Estimate}

Finally, combining Lemmas \ref{chapelle}, \ref{lemma:error}, and \ref{lemma:app}, we 
obtain the following energy norm error estimate.
\begin{thm}\label{apriori} With $(u,\bftheta)$ the exact solution defined by 
(\ref{exactsol}), $(u^h,\bftheta^h)$ the discrete solution defined by (\ref{niteq}), and 
the continuous and discrete shear stress, $\bfzet$ and $\bfzet^h$, defined by 
(\ref{zeta}), we have the estimate
\label{eq:uniform-convergence} 
\[\
\tnorm{\bftheta-\bftheta^h}+
t\norm{\bfzet-\bfzet^h }{L_2(\Omega)}
\leq C h\Bigl(\norm{g}{H^{-1}(\Omega)}+t\norm{g}{L_2(\Omega)}\Bigr) 
\]
uniformly in $t$.
\end{thm}

\section{Numerical examples}

\subsection{Practical implementation}

We focus on a bilinear approximation of the geometry, $\bfx(\widehat{\bfx}) = \bfx_i\widehat{\psi}_i(\widehat{\bfx})$ 
where $\bfx_i$ are the corner node coordinates and
\[
\widehat{\bfpsi}=[(1-\widehat{x})(1-\widehat{y}),\widehat{x}(1-\widehat{y}),\widehat{x}\widehat{y},(1-\widehat{x})\widehat{y}], 
\]
so that
\[
D\bfF_K = \left[\begin{array}{>{\displaystyle}c>{\displaystyle}c}\sum_{i=1}^4 x_i\frac{\partial \widehat{\psi}_i}{\partial \widehat{x}} & \sum_{i=1}^4 x_i\frac{\partial  \widehat{\psi}_i}{\partial \widehat{y}}\\[3mm]
\sum_{i=1}^4 y_i\frac{\partial\widehat{\psi}_i}{\partial \widehat{x}} & \sum_{i=1}^4 y_i\frac{\partial\widehat{\psi}_i}{\partial \widehat{y}}\end{array}\right] .
\]
Then, the covariant map of the rotations is given by
\[
\bftheta (\bfx(\widehat{\bfx})) = D\bfF_K^{-\text{T}}\widehat{\bftheta}(\widehat{\bfx}) ,
\]
or, inversely,
\begin{equation}\label{nummap}
\widehat{\bftheta}(\widehat{\bfx}) = D\bfF_K^{\text{T}}\bftheta(\bfx(\widehat{\bfx}))   .
\end{equation}
Computing the parametric derivatives of $\bftheta$ follows from  applying the derivatives to (\ref{nummap}):
\[
\frac{\partial\bftheta}{\partial \widehat{x}}=D\bfF_K^{-\text{T}}\left(\frac{\partial\widehat{\bftheta}}{\partial \widehat{x}}-\left(\frac{\partial}{\partial\widehat{x}}D\bfF_K^{\text{T}}\right)\bftheta \right);
\]
\[
\frac{\partial\bftheta}{\partial \widehat{y}}=D\bfF_K^{-\text{T}}\left(\frac{\partial\widehat{\bftheta}}{\partial \widehat{y}}-\left(\frac{\partial}{\partial\widehat{y}}D\bfF_K^{\text{T}}\right)\bftheta \right),
\]
and the gradient operator in physical coordinates applied to $\bftheta$ is finally computed via
\[
\left[\begin{array}{>{\displaystyle}c} \frac{\partial }{\partial x}\\[3mm] \frac{\partial }{\partial y}\end{array}\right] \bftheta = \left(D\bfF_K^{-\text{T}}\left[\begin{array}{>{\displaystyle}c} \frac{\partial }{\partial \widehat{x}}\\[3mm] \frac{\partial }{\partial \widehat{y}}\end{array}\right]\right)\bftheta .
\]

\subsection{Convergence}

\begin{figure}[h]
\begin{center}
\includegraphics[width=3in]{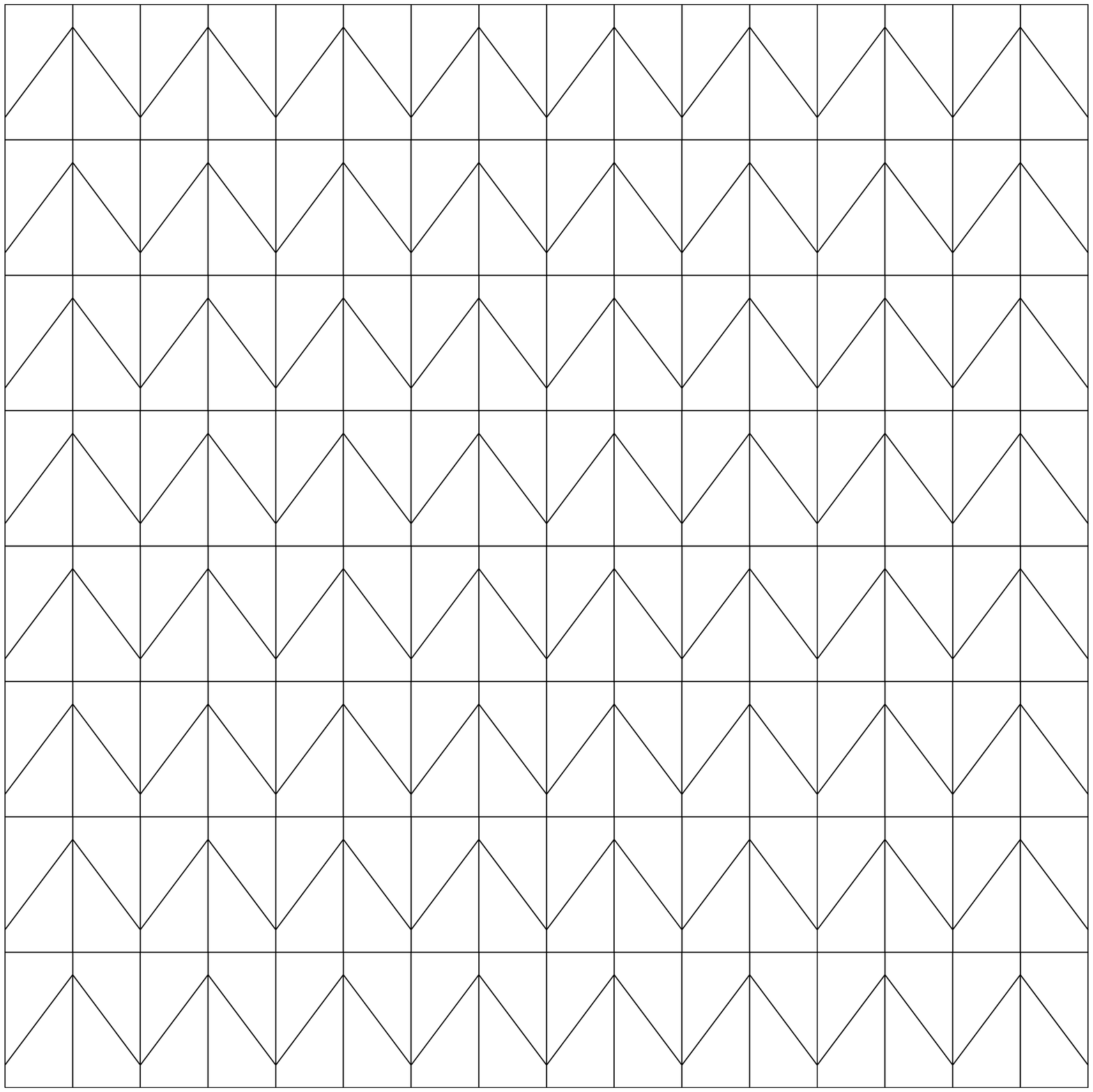}\includegraphics[width=2.8in]{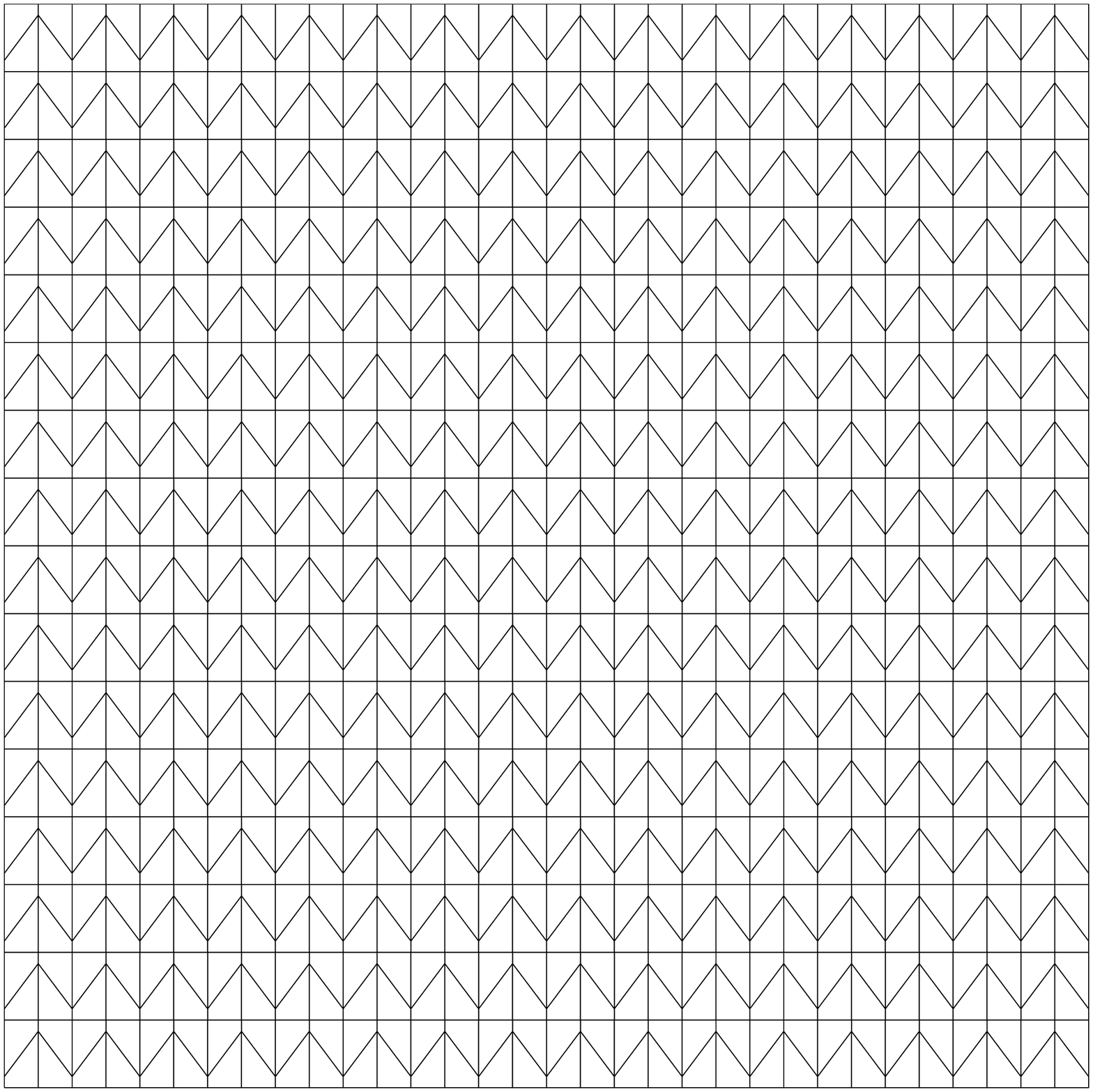}
\end{center}
\caption{First computational mesh and first refined mesh in a sequence.}\label{fig:conv}
\end{figure}

We consider the exact solution to a clamped Reissner--Mindlin plate on the unit square presented
by Chinosi and Lovadina \cite{ChLo95}. They suggested a right-hand side
\begin{align}
g = & {} \frac{E}{12(1-\nu^2)}(12y(y-1)(5x^2-5x+1)
(2y^2(y-1)^2+x(x-1)(5y^2-5y+1))+ \nonumber\\
& 12x(x-1)(5y^2-5y+1)(2x^2(x-1)^2+y(y-1)(5x^2-5x+1))) ,\nonumber
\end{align}
leading to $u = u_0 + u_r$, where
\[
u_0=\frac13 x^3(x-1)^3y^3(y-1)^3 
\]
corresponds to the Kirchhoff solution as $t\rightarrow 0$,
and
\[
u_r=\frac{2t^2}{5(1-\nu)}(y^3(y-1)^3x(x-1)(5x^2-5x+1)+x^3(x-1)^3y(y-1)(5y^2-5y+1)),
\]
and rotations
\[
\theta_x=(y^3(y-1)^3x^2(x-1)^2(2x-1)), \quad \theta_y=(x^3(x-1)^3y^2(y-1)^2(2y-1)).
\]
We let $E=180$ GPa and $\nu=0.3$ and do a study of convergence on a sequence of self-similar trapezoids (following \cite{ArBo02}), as indicated in Fig. \ref{fig:conv}. We consider a continuous $Q_2$--approximation of the displacements, and the rotations in a reference coordinate system, using the covariant map, are, element wise,
\[
\widehat{\theta}^h_{\widehat{x}}\vert_{\widehat{K}} \in \text{span}\{1,\widehat{x},\widehat{y},\widehat{x}\widehat{y},\widehat{y}^2,\widehat{x}\widehat{y}^2\},\quad \widehat{\theta}_{\widehat{y}}^h\vert_{\widehat{K}} \in \text{span}\{1,\widehat{x},\widehat{y},\widehat{x}\widehat{y},\widehat{x}^2,\widehat{y}\widehat{x}^2\} .
\]
For a standard bilinear map the components of ${\widehat{\bftheta}}$ are instead given in physical coordinates. In our implementation, we have
used the same approximating polynomials; thus,
for the bilinear map 
\[
\widehat{\theta}_{{x}}^h\vert_{\widehat{K}} \in \text{span}\{1,\widehat{x},\widehat{y},\widehat{x}\widehat{y},\widehat{y}^2,\widehat{x}\widehat{y}^2\},\quad \widehat{\theta}_{{y}}^h\vert_{\widehat{K}} \in \text{span}\{1,\widehat{x},\widehat{y},\widehat{x}\widehat{y},\widehat{x}^2,\widehat{y}\widehat{x}^2\} .
\]

\begin{figure}[h]
\begin{center}
\includegraphics[width=2in]{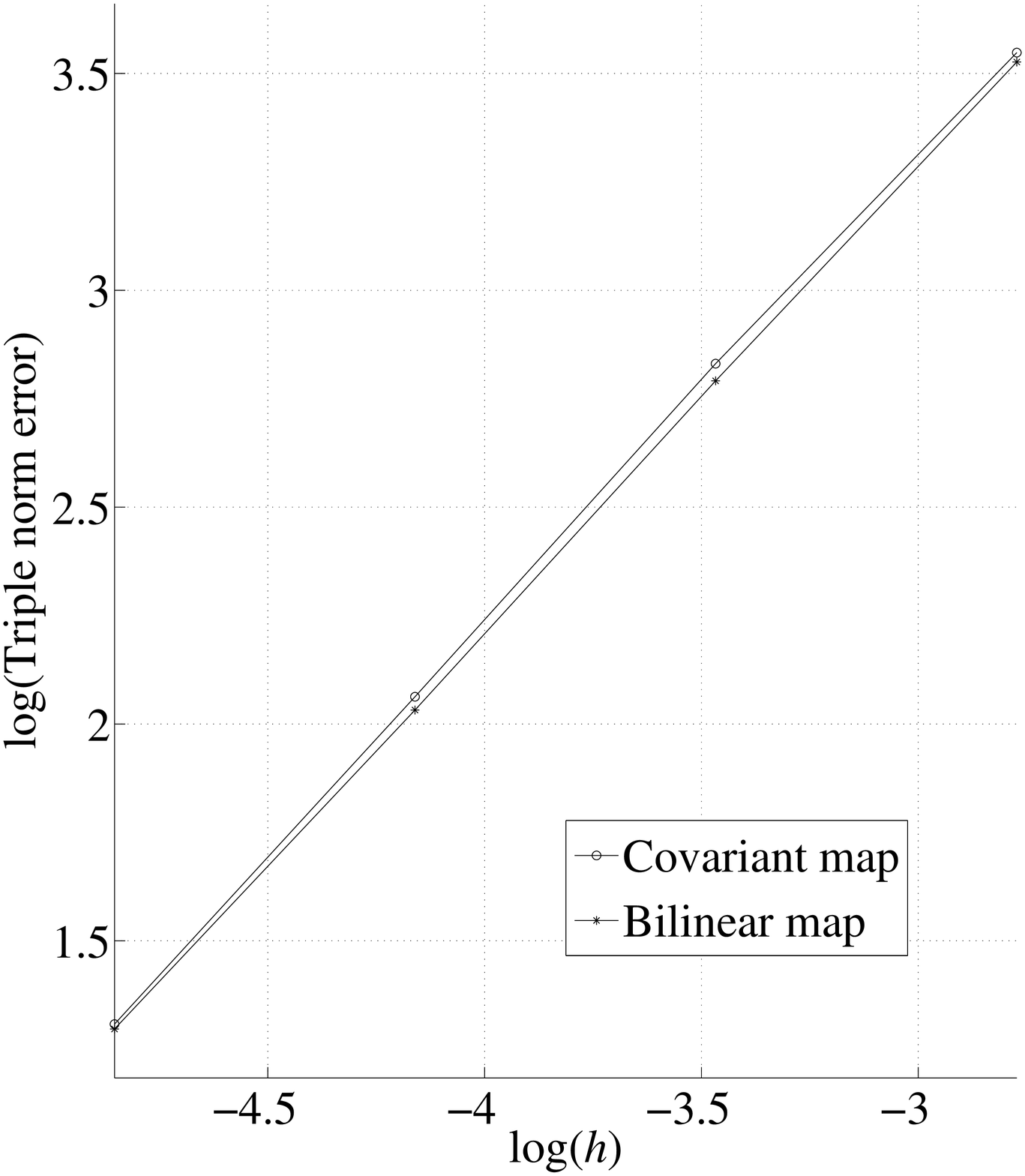}\includegraphics[width=2.5in]{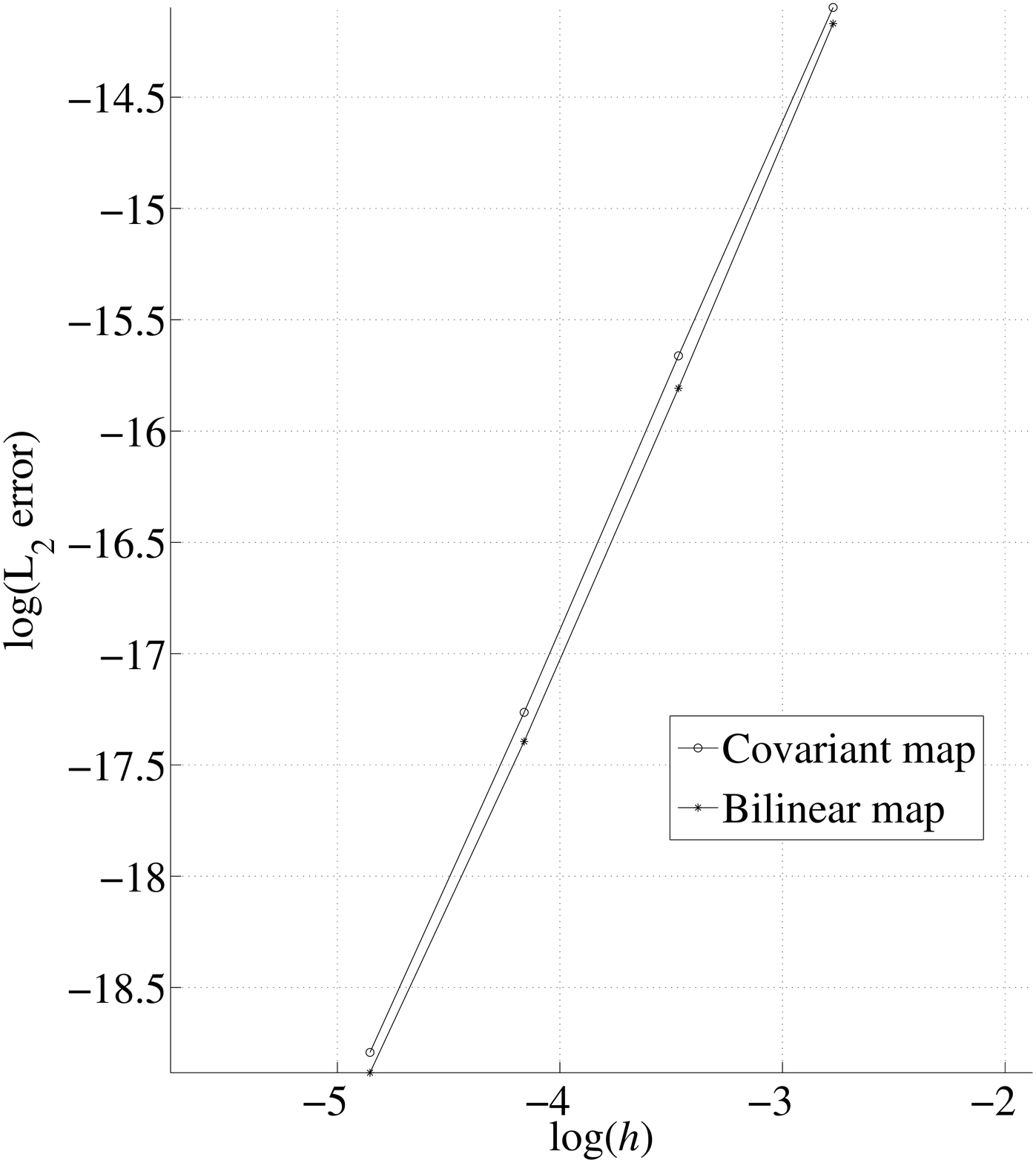}
\end{center}
\caption{Error in triple norm and $L_2(\Omega)$ for $t=10^{-2}$.}\label{fig:conv1}
\end{figure}

\begin{figure}[h]
\begin{center}
\includegraphics[width=2.5in]{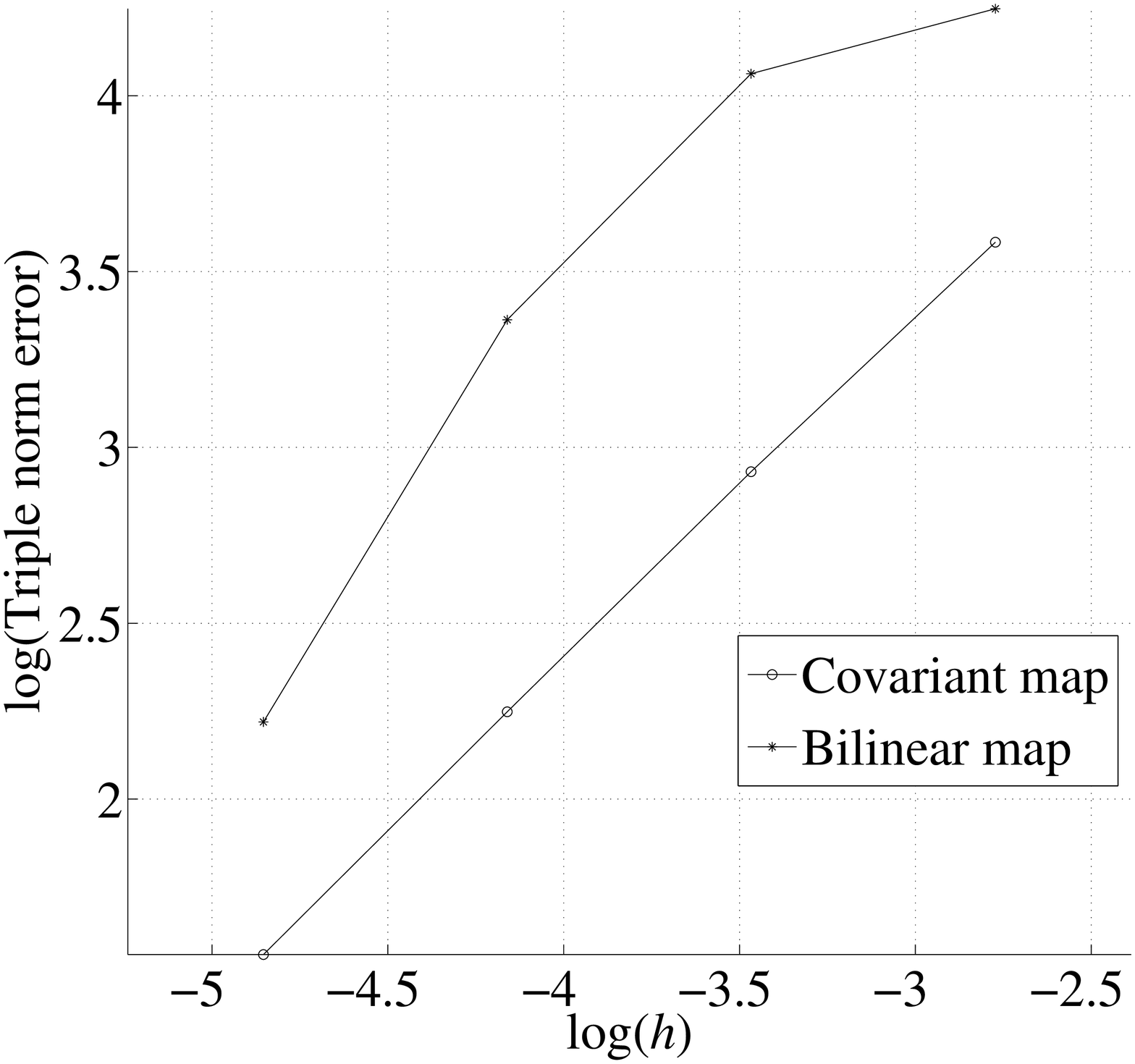}\includegraphics[width=2.5in]{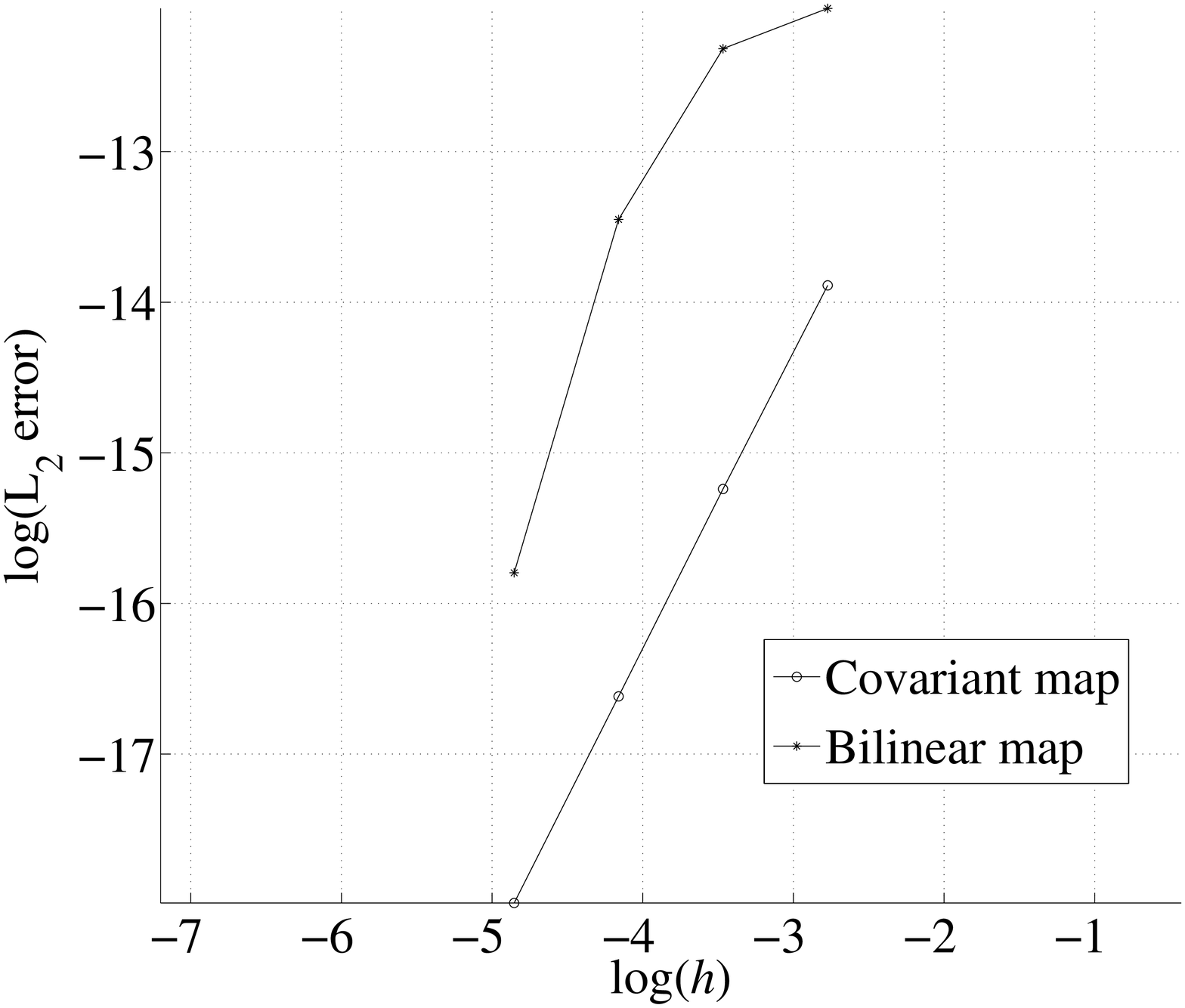}
\end{center}
\caption{Error in triple norm and $L_2(\Omega)$ for $t=10^{-4}$.}\label{fig:conv2}
\end{figure}

The convergence is given for $\tnorm{\bftheta-\bftheta^h}$ and for $\| u-u^h\|_{L_2(\Omega)}$. For $t=10^{-2}$ we observe first and second order convergence, respectively, using
both a standard bilinear map of $\widehat\bftheta$ and the covariant map, cf. Fig. \ref{fig:conv1}. As $t$ becomes smaller, the bilinear map eventually suffers from locking, as illustrated in Fig. \ref{fig:conv2} in the case $t=10^{-4}$. The covariant map is unaffected by the size of $t$. 

\subsection{Locking on a fixed mesh}

\begin{figure}[h]
\begin{center}
\includegraphics[width=3in]{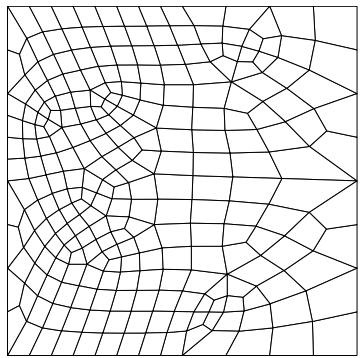}
\end{center}
\caption{Mesh used to illustrate locking.}\label{fig:lock}
\end{figure}

\begin{figure}[h]
\begin{center}
\includegraphics[width=3in]{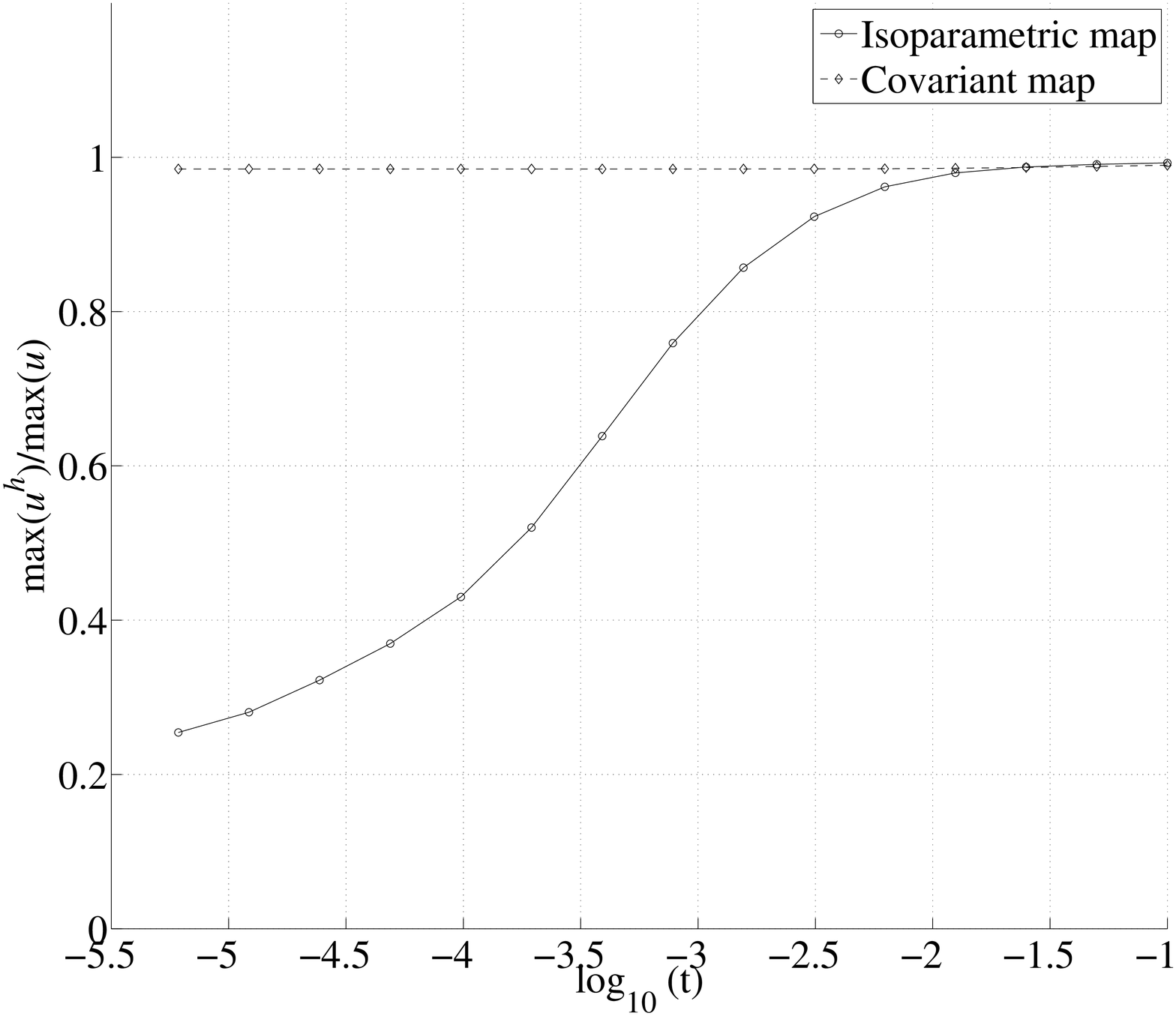}
\end{center}
\caption{Locking using the isoparametric map.}\label{fig:lock1}
\end{figure}

We illustrate the locking problem of a bilinear map further using the unstructured fixed mesh of Fig. \ref{fig:lock}. Using the same problem, approximation, and data as in the previous Section, we plot the ratio of maximum computed displacement to maximum exact displacement (measured in all of the nodes of the $Q_2$ mesh). In Fig. \ref{fig:lock1} we note
the distinct locking effect of using a bilinear map whereas the covariant map is unaffected by the thickness.

\bibliographystyle{plain}
\bibliography{plate}

\newpage



\end{document}